\documentclass{amsart}
\usepackage[all]{xy}

\newtheorem{theorem}{Theorem}[section]

\newtheorem{corollary}[theorem]{Corollary}
\newtheorem{definition}[theorem]{Definition}

\theoremstyle{remark}
\newtheorem{remark}[theorem]{Remark}

\newcommand{\R}{\mathbb{R}}
\newcommand{\C}{\mathbb{C}}

\newcommand{\cO}{\mathcal{O}}
\renewcommand{\a}{\alpha}
\renewcommand{\b}{\beta}

\numberwithin{equation}{section}

\begin{document}

\title{Codimension Theorems for Complete Toric
Varieties}

\author{David Cox}

\address{Department of Mathematics and Computer
Science, Amherst
College, Amherst, MA 01002-5000}

\email{dac@cs.amherst.edu}

\thanks{D.C. would like to thank the Mathematics
Department of the University of Buenos Aires for their
hospitality
during his visits there in 2001 and 2003.}

\author{Alicia Dickenstein}
\address{Departamento de Matem\'atica, F.C.E.\ y N.,
Universidad de
Buenos Aires, Cuidad Universitaria--Pabell\'on I, 1428
Buenos
Aires, Argentina} \email{alidick@dm.uba.ar}
\thanks{A.D. was supported by ANPCYT 03-06568, UBACYT
X-052 and
Conicet, Argentina.}



\begin{abstract}
Let $X$ be a complete toric variety with homogeneous coordinate ring
$S$. In this article, we compute upper and lower bounds for the
codimension in the critical degree of ideals of $S$ generated by
$\dim(X)+1$ homogeneous polynomials that don't vanish simultaneously
on $X$.
\end{abstract}

\maketitle

\section*{Introduction}

In this paper, $X$ will denote a complete toric variety of dimension
$n$.  We will work over $\C$, so that the torus of $X$ is
$(\C^*)^n$. The dual lattices will be denoted $M$ and $N$ as usual,
and the minimal edge generators of the fan of $X$ will be denoted
$n_1,\dots,n_r \in N$.  Corresponding to each $n_i$ we have the
irreducible torus-invariant divisor $D_i$ and the variable $x_i$ in
the homogeneous coordinate ring $S = \C[x_1,\dots,x_r]$, which is
graded by the Chow group $A_{n-1}(X)$.

Consider a torus-invariant Cartier divisor $D = \sum_i a_i D_i$ such
that $\cO_X(D)$ is generated by global sections.  In
Section~\ref{vanishingsec}, we recall a vanishing theorem of Batyrev
and Borisov \cite{bb} that describes $H^i(X,\cO_X(-D))$ in terms of
the polytope
\[
\Delta_D = \{m \in M_\R \mid \langle m,n_i\rangle \ge
-a_i\}.
\]
By Serre duality, we get a description of
\[
H^i(X,\cO_X(D+K)),
\]
where $K$ is the canonical class.

In Section~\ref{codimsec}, we consider homogeneous polynomials
$f_{0},\dots,f_{n} \in S$ of degrees $\a_{0},\dots,\a_{n} \in
A_{n-1}(X)$ which satisfy the following two properties:\ first, the
$f_{i}$ don't vanish simultaneously on $X$, and second, each $\a_i$
lies in $\mathrm{Pic}(X)$ and the corresponding line bundle
$\cO_{X}(\a_{i})$ is generated by global sections.  By abuse of
terminology, we will say that $\a_{i}$ is \emph{globally generated} in
this situation.

Following \cite{ccd}, we define the \emph{critical degree} to be
\begin{equation}
\label{critdeg} \rho = \sum_{i=0}^{n} \a_{i} -
\beta_{0}
\end{equation}
where $\beta_{0} = \sum_{j=1}^{r} \deg(x_{j})$, and we
let
\[
I = \langle f_{0},\dots,f_{n}\rangle \subset S
\]
be the homogeneous ideal generated by the $f_i$.  Then the
\emph{Codimension Question} asks what is the dimension of
$(S/I)_{\rho}$, i.e., what is the codimension of $I_\rho$ in
$S_{\rho}$?  

The main result of Section~\ref{codimsec} computes upper and lower
bounds for this codimension, when the associated family of polytopes
is \emph{essential}. In particular, it is always nonzero.  The
assumption of having an essential family is natural since it is
necessary for the existence of a nontrivial sparse resultant whose
vanishing is equivalent to the condition that the $f_i$ do not have
common zeros on $X$ (see \cite{bernd}). We also give further geometric
conditions on the polytopes under which the codimension at the
critical degree attains the upper bound that we present.

When the divisors involved are big and nef, our results imply that the
codimension is $1$.  The first toric codimension one theorem, proved
in \cite{tr} using ideas from \cite{batyrev}, assumed that the $\a_i$
were all equal to the class of a single ample divisor.  This was
extended in \cite{cd} to the case when the $\a_i$ were positive
integer multiples of a single ample class.  In the general case when
the $\a_i$ are all ample, a codimension one theorem was proved in
\cite{ccd} provided that $X$ was simplicial.  All of these results now
follow from part (1) of Corollary~\ref{bignef}.

The codimension in the critical degree is important because of its
relation to the theory of toric residues \cite{bm, ccd, cd, tr, dk}.
Toric residues are rational functions of the coefficients of the given
polynomials $f_0, \dots, f_n$; moreover, they are rational
hypergeometric functions determined by the lattice points of the
associated family of polytopes, with poles at the resultant locus
\cite{cds}.  As described in \cite{tr}, a homogeneous polynomial $H$
of critical degree gives rise to a rational $n$-form on $X$.  Since
the $f_i$ don't vanish simultaneously on $X$, this $n$-form represents
an element of $H^n(X,\widehat{\Omega}_X^n) \simeq
H^n(X,{\cO}_X(-\beta_0))$.  The toric residue of $H$ is defined to be
the trace of this cohomology class.  If the codimension is one and we
have an explicit element $J$ of critical degree with known residue,
then the computation of the toric residue of $H$ is reduced to writing
$H = c J$ modulo $I$.  When all degrees are equal and ample, a choice
of $J$ is the toric Jacobian \cite{tr} associated to $f_0, \dots,
f_n$.  In case all degrees are ample, explicit elements with residue
equal to $\pm 1$ are known \cite{ccd,dk}, but it is still an open
problem to find such elements in the big and nef but not ample case.

\section{The Vanishing Theorem}
\label{vanishingsec}

Now let $D$ be be torus-invariant Cartier divisor on a
complete
toric variety $X$.  We define the polytope $\Delta_D$
as in the
introduction and we use $\mathrm{int}(\Delta_D)$ to
denote the
relative interior of $\Delta_D$.  Here is the
vanishing theorem
from \cite[Thm.\ 2.5]{bb}.

\begin{theorem}
\label{vanishing} Let $D$ be a torus-invariant Cartier
divisor on
a complete toric variety $X$ and let $\Delta_D$ be the
polytope
defined above. If $\cO_X(D)$ is generated by global
sections,
then:
\begin{enumerate}
\item $H^i(X,\cO_X(-D)) = 0$ for $i \ne
\dim(\Delta_D)$. \item
There is an isomorphism
\[
H^{\dim(\Delta_D)}(X,\cO_X(-D)) \simeq \bigoplus_{m
\in M\cap
\mathrm{int}(\Delta_D)} \C\cdot \chi^{-m}
\]
which is equivariant with respect to the natural
$(\C^*)^n$ action
on each side.
\end{enumerate}
\end{theorem}

\begin{remark} We should note that for a complete
toric variety $X$, the sheaf $\cO_X(D)$ is generated
by its global
sections if and only if $D$ is a nef divisor (this
observation
appears in \cite{mav} and \cite{mus}). Also, in the
terminology of
\cite{mav}, such a divisor $D$ is called
$\ell$-semiample, where
$\ell = \dim(\Delta_D)$.
\end{remark}

Using Serre Duality, we get the following corollary of
Theorem~\ref{vanishing}.  As usual, $K = K_X$ denotes the canonical
divisor of $X$.

\begin{corollary}
\label{cor} Under the same hypotheses as
Theorem~\ref{vanishing},
we have:
\begin{enumerate}
\item $H^i(X,\cO_X(D+K)) = 0$ for $i \ne
n-\dim(\Delta_D)$. \item
There is an isomorphism
\[
H^{n-\dim(\Delta_D)}(X,\cO_X(D+K)) \simeq \bigoplus_{m
\in M\cap
\mathrm{int}(\Delta_D)} \C\cdot \chi^m
\]
which is equivariant with respect to the natural
$(\C^*)^n$ action
on each side.
\end{enumerate}
\end{corollary}

\begin{proof}
We know that $K$ is given by the Weil divisor
$-\sum_{j=1}^{r}
D_{i}$. By Serre Duality,
\[
H^{i}(X,\cO_{X}(D)\otimes_{\cO_{X}}\!\cO_{X}(K))
\simeq
H^{n-i}(X,\cO_{X}(-D))^{*}.
\]
However, given any Weil divisors $E$ and $F$ on $X$,
the natural
map
\[
\cO_{X}(E) \otimes_{\cO_{X}}\!\cO(F) \to
{\cO}_{X}(E+F)
\]
is easily seen to be an isomorphism when $E$ or $F$ is
Cartier.
Since $D$ is Cartier by assumption, the above duality
may be
written
\[
H^{i}(X,\cO_{X}(D+K)) \simeq
H^{n-i}(X,\cO_{X}(-D))^{*}.
\]
By functorality, this is compatible with the
$(\C^{*})^{n}$ action
on everything.  {}From here, the corollary follows
immediately
from Theorem~\ref{vanishing}.
\end{proof}

\begin{remark} When $D$ is big and nef, we have
$\dim(\Delta_D) = n$.
Then Corollary~\ref{cor} implies that
\[
H^i(X,\cO_X(D+K)) = 0, \quad i > 0.
\]
This is the Kawamata-Viehweg vanishing theorem from
\cite{mus}.
\end{remark}

\section{Codimension in the Critical Degree}
\label{codimsec}

As in the introduction, fix $\a_{0},\dots,\a_{n} \in
\mathrm{Pic}(X) \subset A_{n-1}(X)$ such that each
$\a_{i}$ is
globally generated. This implies that each $\a_{i}$
determines a
lattice polytope $\Delta_{i}$ which is well-defined up
to
translation by an element of $M$.

Given polynomials $f_{i} \in S_{\a_{i}}$ (so that
$\deg(f_{i}) =
\a_{i})$, we obtain the homogeneous ideal $I = \langle
f_{0},\dots,f_{n}\rangle \subset S$.  We want to study
the
codimension of $I_\rho$ in $S_\rho$, where $\rho =
\sum_{i=0}^n
\a_i - \beta_0$ is the critical degree
\eqref{critdeg}.

Before stating our main result we need a definition.

\begin{definition} \label{def:essential}
Let $\Delta_{i}$, $i \in I$, be
 polytopes in $\R^{n}$.  Given $J  \subset I$, set
\[
\Delta_{J} = {\textstyle \sum_{j \in J} \Delta_{j}}.
\]
The family $\{\Delta_{i}, i \in I\}$  is called
{\bfseries
essential} if for every $J \subset I$ with $|J| \le
n$, we have
\[
\dim(\Delta_{J}) \ge |J|.
\]
\end{definition}

We will also use the standard notation
\begin{equation}
\label{lstar} l^*(\Delta) = \#(M \cap
\mathrm{int}(\Delta))
\end{equation}
to denote the number of lattice points of a polytope
$\Delta$
which lie in the relative interior of $\Delta$.

Here is our codimension theorem.

\begin{theorem}
\label{codim} Suppose that $X$ is a complete toric
variety of
dimension $n$.  Let $\a_i \in \mathrm{Pic}(X)$, $0 \le
i \le n$,
be globally generated and assume that the
corresponding family of
polytopes $\{\Delta_{i}, 0\le i \le n\}$ is essential.
If $f_{i}
\in S_{\a_{i}}$, $0 \le i \le n$, don't vanish
simultaneously on
$X$, then{\rm :}
\begin{enumerate}
\item The codimension of $I_\rho$ in $S_\rho$
satisfies the
  inequalities
\[
1 + \sum_{\dim(\Delta_\ell)=1} l^*(\Delta_\ell) \le
\dim((S/I)_\rho) \le 1 + \sum_{k=1}^{n-1}
\sum_{\dim(\Delta_J)
=|J| = k} l^*( \Delta_J ),
\]
where in the right-most sum we always assume that $J
\subset
\{0,\dots,n\}$. \item The codimension of $I_\rho$ in
$S_\rho$ is
given by
\[
\dim((S/I)_\rho) = 1 + \sum_{k=1}^{n-1}
\sum_{\dim(\Delta_J) =|J|
= k} l^*(\Delta_J)
\]
if one of the following is satisfied for every $J$
with $1 \le |J|
\le n-2${\rm :}
\begin{enumerate}
\item $\dim(\Delta_J) \ne |J|+1$. \item
$\dim(\Delta_J) = |J|+1$
but $\Delta_J$ has no interior
 lattice points.
\item $\dim(\Delta_J) = |J|+1$ but $\dim(\Delta_{J\cup
I}) > |J| +
  |I|$ for nonempty subsets $I \subset \{0,\dots,n\}$
such that $I
  \cap J = \emptyset$ and $|J| + |I| < n$.
\end{enumerate}
\end{enumerate}
\end{theorem}

\begin{proof}
Since $f_{i}$ is a global section of
$\cO_{X}(\a_{i})$, we get a
Koszul complex
\[
0 \to \cO_{X}(-{\textstyle\sum_i} \a_i) \to \dots \to
\bigoplus_{i} \cO_{X}(-\a_{i}) \to {\cO}_{X} \to 0
\]
which is exact since each $\cO_{X}(\a_{i})$ is locally
free and
the $f_{i}$ have no common zeros.

By hypothesis, each sheaf in the Koszul complex is
locally free.
Since the sheaf $\underline{\mathrm{Tor}}_i^{\cO_X}
(\mathcal{E},\mathcal{F})$ vanishes whenever $i > 0$
and
$\mathcal{E}$ or $\mathcal{F}$ is locally free, it
follows that
the Koszul sequence remains exact after tensoring with
$\cO(\rho)
= \cO(\a_{0}+ \dots +\a_{n} -\beta_{0})$ (which need
not be
locally free).  This gives the exact sequence
\[
0 \to \cO_{X}(-{\textstyle\sum_i}
\a_i)\otimes_{\cO_{X}}\!\cO_X(\rho) \to \cdots \to
\bigoplus_{i}
\cO_{X}(-\a_{i}) \otimes_{\cO_{X}}\!\cO_X(\rho) \to
\cO_X(\rho)
\to 0.
\]
Since the $\a_{i}$ all come from Cartier divisors, the
reasoning
used in the proof of Corollary~\ref{cor} implies that
we can write
this exact sequence as
\[
0 \to \cO_{X}(-\beta_{0}) \to \bigoplus_{i}
\cO_{X}(\a_{i}-\beta_{0}) \to \dots \to
{\cO}_{X}(\a_{0}+\dots+\a_{n}-\beta_{0}) \to 0.
\]
Using $\rho = \a_{0}+\dots+\a_{n}-\beta_{0}$, this
becomes
\begin{equation}
\label{imptex} 0 \to
\underbrace{\cO_{X}(-\beta_{0})_{\strut}}_{\mathcal{F}^0}\!
\to
\underbrace{\bigoplus_{i}
\cO_{X}(\a_{i}-\beta_{0})}_{\mathcal{F}^1} \to\!
\cdots\! \to
\underbrace{\bigoplus_{k} \cO_{X}(\rho -
\a_{k})}_{\mathcal{F}^n}
\to
\underbrace{\cO_{X}(\rho)_{\strut}}_{\mathcal{F}^{n+1}}\!
\to
0.
\end{equation}
We will now study the hypercohomology of the complex
$\mathcal{F}^\bullet$.

Since \eqref{imptex} is exact, the hypercohomology
$\mathbb{H}^*(X,\mathcal{F}^\bullet)$ vanishes
identically.  Hence
we get a spectral sequence
\[
E_1^{p,q} = H^q(X,\mathcal{F}^p) \Rightarrow 0
\]
where the differential $d_1^{p,q} : E_1^{p,q} \to
E_1^{p+1,q}$ is
induced by the map $\mathcal{F}^p \to
\mathcal{F}^{p+1}$ of the
Koszul complex \eqref{imptex}.

We compute $E_1^{p,q}$ as follows.  When $p = 0$,
\begin{equation}
\label{e0i} \dim(E_1^{0,q}) =
\dim\big(H^{q}(X,\cO_{X}(-\beta_{0})\big) =
\begin{cases} 0 & q < n \\ 1 & q = n, \end{cases}
\end{equation}
where $q=n$ uses the isomorphism
$H^{n}(X,\cO_{X}(-\beta_{0})) =
H^{n}(X,\cO_{X}(K)) \simeq \C$ given by the trace map.

When $p > 0$, we have
\[
\mathcal{F}^p = \bigoplus_{|J| = p} \cO_X(\a_J -
\beta_0),\quad
\a_J = {\textstyle{\sum_{j\in J}}} \a_j.
\]
Then Corollary~\ref{cor} and \eqref{lstar} imply that
\begin{equation}
\label{e1} \dim(E_1^{p,q}) = \sum_{|J| = p}
\dim\big(H^q(X,\cO_X(\a_J - \beta_0))\big) = \sum_{|J|
= p,
\dim(\Delta_J) = n-q} l^*(\Delta_J).
\end{equation}
Our assumption that $\{\Delta_i, 0 \le i \le n\}$ is
essential
implies that $\dim(\Delta_J) \ge |J|$ when $|J| \le
n$.  Combining
this with $|J| = p$, $\dim(\Delta_J) = n-q$ from
\eqref{e1}, we
see that $n-q \ge p$ when $p \le n$ in the last
summation of
\eqref{e1}.  Thus
\begin{equation}
\label{ssvanish} E_1^{p,q} = 0\ \text{when}\ p+q > n,\
\text{except for}\ E_1^{n+1,0} = S_\rho.
\end{equation}
When $n = 4$, \eqref{e0i} and \eqref{ssvanish} give the picture shown
in {\scshape Figure}~1 on the next page of all possible nonzero
differentials in the spectral sequence

\begin{figure}[t]
\caption{Possible nonzero differentials when $n = 4$}
\[
\xymatrix{
\bullet  \ar[ddddrrrrr]|{d_5} \\
\circ & \bullet  \ar[dddrrrr]|{d_4} \\
\circ & \bullet \ar[r]|{d_1} \ar[drr]|{d_2}
\ar[ddrrr]|{d_3} & \bullet \ar[ddrrr]|{d_3} \\
\circ & \bullet \ar[r]|{d_1} \ar[drr]|{d_2} & \bullet
\ar[drr]|{d_2}
\ar[r]|{d_1} & \bullet \ar[drr]|{d_2}\\
\circ & \bullet \ar[r]|{d_1} & \bullet \ar[r]|{d_1} &
\bullet
\ar[r]|{d_1} & \bullet \ar[r]|{d_1} & \bullet }
\]
\end{figure}
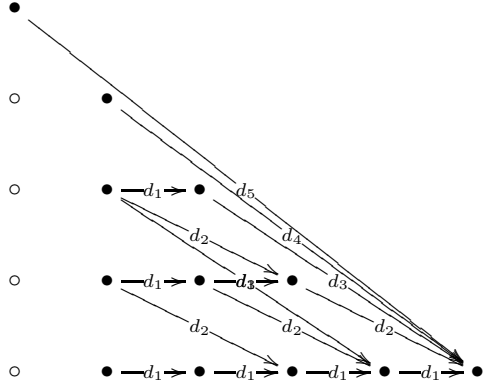

Now consider $E_r^{n+1,0}$.  By \eqref{imptex},
$E_1^{n,0} \to
E_1^{n+1,0} \to E_1^{n+2,0}$ is
\[
\bigoplus_k S_{\rho-\a_k} \longrightarrow S_\rho
\longrightarrow
0.
\]
Since this comes from the Koszul complex of the $f_i$,
the image
of the first map is $I_\rho$.  It follows that
\[
E_2^{n+1,0} = (S/I)_\rho.
\]
All differentials starting from $E_2^{n+1,0}$ obviously vanish, and
the only differentials which can map to this position are those on the
``diagonal'' $p+q = n$ as in {\scshape Figure}~1.  Furthermore, on the
diagonal, the only nonzero differentials $d_r^{n-q,q}$ are
$d_{q+1}^{n-q,q} : E_{q+1}^{n-q,q} \to E_{q+1}^{n+1,0}$.  Since the
spectral sequence converges to zero, it follows easily that
\[
\dim((S/I)_\rho) \le  \sum_{q=1}^{n}
\dim(E_{1}^{n-q,q}).
\]
By \eqref{e0i} and \eqref{e1}, we get the upper bound
of part (1)
of the theorem.  As for the lower bound, note that the
differentials
\[
d_{n+1}^{0,n} : E_{n+1}^{0,n} \to E_{q+1}^{n+1,0}\quad
\text{and}
\quad d_{n}^{1,n-1} : E_{n}^{1,n-1} \to
E_{q+1}^{n+1,0}
\]
must be injective since the spectral sequence
converges to $0$ and
nothing can map to these positions.  (In {\scshape
Figure}~1,
these correspond to the differentials $d_5$ and
$d_4$.)  Using
\eqref{e0i} and \eqref{e1} again, we get the lower
bound of part
(1).

Turning to part (2) of the theorem, suppose that every
differential $d_r$ mapping to $E_r^{p,q}$ is zero for
$p+q = n$
and $q > 0$.  Since the spectral sequence converges to
zero, this
implies that
\[
\dim((S/I)_\rho) =  \sum_{q=1}^{n}
\dim(E_{1}^{n-q,q}).
\]
{}From this, we easily get the desired formula for the
codimension. Hence it suffices to prove that these
differentials
vanish when condition (a), (b) or (c) is satisfied by
any $J$ with
$1 \le |J| \le n-2$.  A differential mapping to the
diagonal $p+q
= n$ originates from the ``sub-diagonal'' $p+q = n-1$.
 Thus we
need to prove that all differentials
\[
d_r^{p,q} : E_r^{p,q} \longrightarrow E_r^{p+r,q-r+1}
\]
vanish when $p+q = n-1$ and $q-r+1 > 0$ (i.e., $r \le
q$).

To analyze this, first note that for arbitrary
$(p,q)$,
Corollary~\ref{cor} implies that
\begin{equation}
\label{E1PQ} E_1^{p,q} = \bigoplus_{|J| =
p,\dim(\Delta_J) = n-q}
H_J,
\end{equation}
where
\begin{equation}
\label{hdescription} H_J = H^q(X,\cO_X(\a_J - \b_0))=
\bigoplus_{m
\in M\cap
  \mathrm{int}(\Delta_J)} \C\cdot \chi^m.
\end{equation}
In particular, when $p+q = n-1$, the description of
$E_1^{p,q}$
becomes
\[
E_1^{p,q} = \bigoplus_{\dim(\Delta_J) = |J|+1 = n-q}
H_J.
\]
The subsets $J \subset \{0,\dots,n\}$ in this direct
sum satisfy:
\begin{equation}
\label{E1PQ1} |J| = p = n-q-1 \le n-2
\quad\text{and}\quad
\dim(\Delta_J) = n-q =
 |J| + 1.
\end{equation}
Every $J$ in the discussion which follows will satisfy
\eqref{E1PQ1}.

Now consider $d_1^{p,q} : E_1^{p,q} \to E_1^{p+1,q}$
for $p+q =
n-1$ and $q \ge 1$.  By \eqref{hdescription} and
\eqref{E1PQ1}, we
see that $H_J = 0$ when $J$ satisfies conditions (a)
or (b) of
part (2) of the theorem.  Now suppose $J$ satisfies
condition (c).
Since $d_1$ comes from the Koszul complex,
\eqref{E1PQ} for
$(p+1,q)$ shows that $d_1^{p,q}$ restricted to $H_J$
is a map
\[
d_1^{p,q} : H_J \longrightarrow \bigoplus_{i \notin J,
  \dim(\Delta_{J\cup\{i\}}) = n-q} H_{J\cup\{i\}}.
\]
By condition (c), we know that
$\dim(\Delta_{J\cup\{i\}}) > |J| +
1$ for $i \notin J$.  But $|J| = n-q-1$ by
\eqref{E1PQ1}, so that
$\dim(\Delta_{J\cup\{i\}}) > n-q$.  Comparing this to
the above
description of $d_1^{p,q}$ shows that $d_1^{p,q} = 0$
when $p+q =
n-1$ and $q \ge 1$, as claimed.

Next consider $d_2^{p,q} : E_2^{p,q} \to
E_2^{p+2,q-1}$, where
$p+q = n-1$ and $q \ge 2$.  To understand this map, we
will recall
its definition which follows from the Snake Lemma. Let
$\mathcal{U}$ be a Leray covering of $X$, so that the
cohomology
groups $H^q(X, \mathcal{F}^p)$ can be identified with
the \v{C}ech
cohomology of $X$ with respect to this covering. As
usual, we call
$\delta$ the maps between cochain groups.  Also let
\[
Z^q(\mathcal{U},\mathcal{F}^p) \subset
C^q(\mathcal{U},\mathcal{F}^p)
\]
denote the corresponding \v{C}ech group and its
subgroup of
cocycles. Then we compute $d_2$ using the following
diagram:
\[
\begin{matrix}
&&&&& C^{q-1}(\mathcal{U},\mathcal{F}^{p+1}) &
\xrightarrow{d_1} &
C^{q-1}(\mathcal{U},\mathcal{F}^{p+2})\\
&&&&& \delta \downarrow \phantom{\delta} &&\\
&&& Z^q(\mathcal{U},\mathcal{F}^p) & \xrightarrow{d_1}
&
  Z^q(\mathcal{U},\mathcal{F}^{p+1}) &&\\
&&& \downarrow && \downarrow &&\\
0 \to & \mathrm{ker}(d_1) & \to & H^q(X,\mathcal{F}^p)
&
\xrightarrow{d_1}& H^q(X,\mathcal{F}^{p+1}) &&\\
& \downarrow && \downarrow && \downarrow &&\\
& E_2^{p,q} && 0 && 0 &&\\
& \downarrow &&&&&&\\
& 0 &&&&
\end{matrix}
\]
In this diagram, $s \in E_2^{p,q}$ lifts to an element
of
$E_1^{p,q} = H^q(X,\mathcal{F}^p)$ represented by $s'
\in
Z^q(\mathcal{U},\mathcal{F}^p)$.  Then
\begin{equation}
\label{ABC} d_1(s') \in
Z^q(\mathcal{U},\mathcal{F}^{p+1})
\end{equation}
represents $0 \in H^q(X,\mathcal{F}^{p+1})$.  Thus
$d_1(s')$ lifts
to $Z^{q-1}(\mathcal{U},\mathcal{F}^{p+1})$, i.e.,
$d_1(s') =
\delta(s'')$ for some $s'' \in
Z^{q-1}(\mathcal{U},\mathcal{F}^{p+1})$.  It is then
easy to see
that
\begin{equation}
\label{BCD} d_1(s'') \in
Z^{q-1}(\mathcal{U},\mathcal{F}^{p+2}),
\end{equation}
so that $d_1(s'')$ represents an element of
$E_1^{p+2,q-1} =
H^{q-1}(X,\mathcal{F}^{p+2})$.  This is killed by
$d_1$ and
represents $d_2^{p,q}(s) \in E_2^{p+2,q-1}$.  It
follows that
$d_2$ is constructed by applying $d_1$ twice to
suitable liftings.

Since $d_1^{p,q}$ vanishes, $E_2^{p,q}$ is a quotient
of
$E_1^{p,q}$. So we may assume that $s \in E_2^{p,q}$
comes from
$H_J$ for some $J$ satisfying \eqref{E1PQ1}.  As
before, this
vanishes if $J$ satisfies conditions (a) or (b) of
part (2).  Now
assume $J$ satisfies condition (c).  Then as above $s$
gives $s'
\in Z^q(\mathcal{U},\mathcal{F}^p)$ where we only use
the summand
of $\mathcal{F}^p$ corresponding to $J$. Since $d_1$
comes from
the Kozsul complex, we see from \eqref{ABC} that
$d_1(s')$ only
involves summands of $\mathcal{F}^{p+1}$ corresponding
to $J \cup
\{i\}$ for $i \notin J$.  Then the same is true for
$s''$, so that
in \eqref{BCD}, $d_1(s'')$ only involves summands of
$\mathcal{F}^{p+1}$ corresponding to $J \cup \{i,j\}$
for
$i,j\notin J$ and $i \ne j$.

This shows that $d_2^{p,q}(s) \in E_2^{p+2,q-1}$ is
represented by
an element of $E_1^{p+2,q-1}$ involving only summands
$H_{J \cup
\{i,j\}}$ for $i,j\notin J$ and $i \ne j$.  Yet to
actually appear
in $E_1^{p+2,q-1}$, the summand must satisfy $n-(q-1)
=
\dim(\Delta_{J\cup \{i,j\}})$ by \eqref{E1PQ}.  This
is impossible
since condition (c) implies that
\[
\dim(\Delta_{J\cup \{i,j\}}) > |J| + 2 = (n-q-1) + 2 =
n-(q-1).
\]
It follows that $d_2^{p,q} = 0$ when $p+q = n-1$ and
$q \ge 2$.

The argument for general $r$ is similar.  Here, $p+q =
n-1$, $q
\ge r$, and $d_1^{p,q}, \dots, d_{r-1}^{p,q}$ vanish. 
Then
$E_r^{p,q}$ is a quotient of $E_1^{p,q}$, so we may
assume that $s
\in E_r^{p,q}$ comes from $H_J$ for $J$ as in
\eqref{E1PQ1}.
Since $d_1$ comes from the Koszul complex and $d_r$ is
obtained by
applying $d_1$ $r$ times to suitable liftings, we see
that
$d_r^{p,q}(s) \in E_r^{p+r,q-r+1}$ is represented by
an element of
$E_1^{p+r,q-r+1}$ involving only summands $H_{J \cup
I}$ for $|I|
= r$ and $I\cap J = \emptyset$.  In order appear in
$E_1^{p+r,q-r+1}$, the summand must satisfy $n-(q-r+1)
=
\dim(\Delta_{J\cup I})$ by \eqref{E1PQ}.  This is
impossible since
condition (c) implies that
\[
\dim(\Delta_{J\cup I}) > |J| + |I| = (n-q-1) + r =
n-(q-r+1).
\]
It follows that $d_r^{p,q} = 0$ when $p+q = n-1$ and
$q \ge r$.
This completes the proof of the theorem.
\end{proof}

Theorem~\ref{codim} has the following corollaries.

\begin{corollary}
\label{interior} Let $\Delta_0,\dots,\Delta_n$ be an
essential
family of lattice polytopes in $\R^n$.  Then the
Minkowski sum
$\Delta_0 + \dots + \Delta_n$ has at least one
interior lattice
point.
\end{corollary}

\begin{proof}
If $X_\Delta$ is the toric variety determined by
$\Delta =
\Delta_0 + \dots + \Delta_n$, then $\dim(S_\rho)$ is
the number of
interior lattice points of $\Delta$.  The lower bound
given by
Theorem~\ref{codim} implies that $\dim(S_\rho) \ge
\dim((S/I)_\rho) \ge 1$, and the corollary follows.
\end{proof}

\begin{remark} Corollary~\ref{interior} may fail if
$\Delta_0,\dots,\Delta_n$ are not essential.  For
example, suppose
that $\Delta_1+\dots+\Delta_{n}$ has dimension $n-1$
and
$\Delta_0$ is an interval of length one relative to an
integral
linear functional constant on the affine hyperplane
containing
$\Delta_1+\dots+\Delta_{n}$.  It is easy to see that
$\Delta_0+\dots+\Delta_{n}$ has no interior lattice
points.  We
are grateful to G\"unter Ziegler for this observation.
\end{remark}

\begin{corollary}
\label{bignef} Suppose that $X$ is a complete toric
variety of
dimension $n$.  Let $\a_i \in \mathrm{Pic}(X)$, $0 \le
i \le n$,
be globally generated and assume that $f_{i} \in
S_{\a_{i}}$, $0
\le i \le n$ don't vanish simultaneously on $X$. 
Then:
\begin{enumerate}
\item If the polytopes $\Delta_{i}$ all have dimension
$n$, then
\[
\dim((S/I)_{\rho}) = 1.
\]
\item If $n = 2$ and $\{\Delta_0,\Delta_1,\Delta_2\}$
is
essential, then
\[
\dim((S/I)_{\rho}) = 1 + \sum_{\dim(\Delta_i) = 1}
l^*(\Delta_i).
\]
\end{enumerate}
\end{corollary}

\begin{proof} This follows immediately from part (2)
of
Theorem~\ref{codim}.
\end{proof}

\begin{remark} In part (1) of Corollary~\ref{bignef}, the hypothesis
  that $\a_i$ is globally generated and $\Delta_i$ has dimension $n$
is equivalent to assuming that $\a_i$ the class of a big and nef
divisor on $X$.
\end{remark}

When $n = 3$, the conditions of part (2) of Theorem~\ref{codim} are
equivalent to the assumption that if $\dim(\Delta_i) = 2$, then either
$\Delta_i$ has no interior lattice points or $\dim(\Delta_i +
\Delta_j) = 3$ for all $j \ne i$.  It follows that Theorem~\ref{codim}
computes the codimension in the critical degree for many but not all
cases of essential supports when $X$ has dimension $3$.

There is one case where further general results are possible.

\begin{theorem}
\label{moregeneral} Suppose that $X$ is a complete
toric variety
of dimension $n \ge 3$. Let $\a_i \in
\mathrm{Pic}(X)$, $0 \le i
\le n$, be globally generated and assume that the
corresponding
family of polytopes $\{\Delta_{i}, 0\le i \le n\}$ is
essential
and satisfies
\begin{equation}
\label{restrictdelta} \dim(\Delta_i) \in \{1
,n-1,n\},\quad i =
0,\dots,n.
\end{equation}
Let $f_{i} \in S_{\a_{i}}$, $0 \le i \le n$ and assume
the
following two conditions:
\begin{enumerate}
\item The $f_i$ don't vanish simultaneously on $X$.
\item For
every $J \subset \{0,\dots,n\}$ with $\dim(\Delta_{J})
= |J| =
n-1$ and $\dim(\Delta_i) = n-1$ for at least one $i
\in J$, the
equations on $X_{\Delta_{J}}$ given by $f_{j} = 0$ for
$j \in J$
have only finitely many solutions, all of which lie in
the torus
of $X_{\Delta_{J}}$.
\end{enumerate}
Then the codimension of $I_\rho$ in $S_\rho$ is given
by
\begin{equation}
\label{genfor} \dim((S/I)_{\rho}) = 1 +
\sum_{\dim(\Delta_J) = |J|
< n}\bigg( \sum_{\mathcal{J} \subset J}
(-1)^{|J|-|\mathcal{J}|}\,
l_{|J|}^*({\textstyle \sum_{\ell \in \mathcal{J}}
P_\ell})\bigg),
\end{equation}
where, for a lattice polytope $\Delta$ and integer $k
\ge 0$,
\[
l^*_k(\Delta) = \begin{cases} 0 & \dim(\Delta) \ne k\\
l^*(\Delta)
&
  \dim(\Delta) = k.\end{cases}
\]
\end{theorem}

The condition \eqref{restrictdelta} is still rather restrictive, yet
examples with $n = 4$ show that the formula \eqref{genfor} can fail
when we omit \eqref{restrictdelta}.  Furthermore, when we do assume
\eqref{restrictdelta}, condition (2) on the $f_i$ is probably
unnecessary, yet we can't figure out how to prove the theorem without
using this hypothesis.  For these reasons, we omit the proof of
Theorem~\ref{moregeneral}.

\end{document}